\documentclass[preprint,authoryear,12pt]{elsarticle}

\usepackage{amssymb}

\journal{Ocean Modelling}
\usepackage{color}
\usepackage{float}
\usepackage{amsmath}
\usepackage{amssymb}
\floatname{algorithm}{Scheme}

\usepackage{caption}
\usepackage{listings}
\usepackage{lineno}

\newtheorem{theorem}{Theorem}[section]

\newenvironment{proof}[1][Proof]{ \begin{trivlist}
\item[\hskip \labelsep {\bfseries #1}]}{\end{trivlist}}

\newenvironment{remark}[1][Remark]{\begin{trivlist}
\item[\hskip \labelsep {\bfseries #1}]}{\end{trivlist}}

\begin{document}
\begin{frontmatter}



\title{\textbf{A Revised  Scheme to Compute  Horizontal  Covariances in an   Oceanographic 3D-VAR  Assimilation System} }

\author[label1]{R. Farina}
\author[label1]{S. Dobricic}
\author[label1]{A. Storto}
\author[label1]{S. Masina}
\author[label2]{and S. Cuomo}
 \address[label1]{Centro Euro-Mediterraneo sui Cambiamenti Climatici,\\ Bologna, Viale Aldo Moro 44, Italy.\\
  Email: raffaele.farina@cmcc.it - Telephone: +39 051 378267 }
 \address[label2]{Department of Mathematics and Applications, University of Naples ''Federico II''.\\
 Via Cinthia, 80126, Napoli, Italy.\\ 
 Email: salvatore.cuomo@unina.it -  Telephone:  +39 081675624}

\begin{abstract}
 {We} propose an improvement  of an oceanographic  three dimensional variational assimilation scheme (3D-VAR), named OceanVar,
by introducing a recursive   filter (RF) with the third order of accuracy  (3rd-RF), instead  of a RF with first order of accuracy (1st-RF), to approximate horizontal Gaussian covariances.
An advantage of the proposed scheme is that the CPU's time  can be substantially reduced with benefits on the  large scale applications. 
Experiments estimating the impact of 3rd-RF  are performed  by assimilating  oceanographic data  in two realistic oceanographic applications. The results evince  benefits in terms of assimilation  process computational time,  accuracy of the Gaussian correlation modeling, and  show   that the 3rd-RF  is a suitable tool for operational data assimilation.

\end{abstract}


\end{frontmatter}


\section{Introduction}
\noindent Ocean data assimilation is a crucial task in operational oceanography, responsible for optimally combining observational measurements and a prior knowledge of the state of the ocean in order to provide initial conditions for the forecast model. How the informative content of the observations is spread horizontally in space depends on the operator used to model horizontal covariances.
 The {three}-dimensional variational data assimilation scheme called {OceanVar \citep{Dobricic} represents horizontal covariances of background errors in temperature and salinity by  approximate  Gaussian functions that depend only on the horizontal distance between the two model points}.  In the framework of the Optimal Interpolation, where the analysis is found by using only the nearest observations, the calculation of the Gaussian function can be made directly from the distances between the model point and the typically small number of nearby observations. This kind of solution may be impractical in the variational framework where it is necessary to calculate the covariances between each pair of model points in the horizontal. Instead, variational schemes often use linear operators that approximate the Gaussian function   \citep[e.g.,][]{Weaver}.\\
In meteorology,   \citet{Lorenc} approximated the Gaussian function by applying one-dimensional recursive filters (RF) with the first-order accuracy successively in the two perpendicular directions.  In oceanography,  \citet{Weaver}   used the explicit integration of the two-dimensional diffusion equation.  \citet{Purser}  developed higher order recursive filters  for use in atmospheric models. The OceanVar scheme described   in  \citet{Dobricic}  used the  RF of the first order  (1st-RF)  with {imaginary sea points  for the processing
on the coast}. \citet{Mirouzea} implemented the implicit integration of one-dimensional diffusion equations. \\
Successive applications of one-dimensional recursive filters or implicit integrations of the diffusion equation in the two perpendicular directions are much more computationally efficient than the explicit integration of the diffusion equation   \citep{Mirouzea}. However, the first order accurate operators used in most of these schemes still require several iterations to approximate the Gaussian function. For example, the  1st-RF in OceanVar  generally  applies 5 iterations and is the computationally most demanding part.\\ 
Recursive filters with higher order accuracy require more operations for each iteration, but  only one iteration might be enough to accurately approximate the Gaussian function. {Generally high order recursive filters can be obtain by means of different strategies  e.g. in  \citet{Purser}, \citet{Young} and \citet{Deriche}. The main difference among them  is the mathematical methodology used to obtain the filter coefficients. In meteorology,  \citet{Purser} resolve an  inverse problem  with   exponential matrix of finite differences operator  approximating the second derivative  $ d^2 /dx^2$  on a line
grid of uniform spacing  $\delta x$. They use truncated  Taylor  expansion   to approximate the exponential matrix  and  obtain the filter coefficients through the  $LL^T$ factorization of the  result  approximation. The  degree  of the Taylor polynomial is  the order of RF obtained. \\ 
In this study we  develop a RF of the third order accuracy (3rd-RF),{still with the use of the imaginary points for treatment on the coast}, that  needs only one iteration to approximate the Gaussian function and  allows different length scales.
 Our approach, based on \citet{Young},  determines the filter coefficients of a 3rd-RF by  the  matrix-vector multiplication of gaussian operator for a input field,  using  a  known  rational approximation of the gaussian function. Note that this strategy has been so far exploited only in signal processing, and represents a completely novel methodology in geophysical data assimilation. Furthermore,  we  compare the 3rd-RF performance with those of the existing 1st-RF on two different configurations  of OceanVar:  Mediterranean Sea and Global Ocean}. The new filter should be at least as accurate as the existing one  and it should execute more rapidly on massively parallel computers. Tests on parallel computer architectures are especially important because higher order accurate filters compute the solution from several nearby points, and as  a consequence, transfer more data among processors eventually becoming  less efficient. \\
Section 2 gives a general description of the existing OceanVar data assimilation scheme. Section 3 demonstrates in detail the development of the  3rd-RF. It also provides all numerical values and the method to calculate the coefficients of the filter  with different length scales. Moreover
we  give an estimate of approximation error between the result of a
RF and the real Gaussian convolution. In Sections 4 and 5   the filter is applied in the operational version of OceanVar used respectively  in the Mediterranean Sea \citep{Pinardis} and Global Ocean \citep{mwr}.  Its performance is compared to the performance of 1st-RF. In Section 5 we present the conclusions and indicate the future directions of the development of the operator for the horizontal covariances.   

\section{General Description}
\subsection{The OceanVar Computational Kernel }
\noindent The computational kernel of the  {OceanVar}
data assimilation scheme  is based on the following regularized constrained least square problem:
\begin{equation}
 \min_D\{J(           x) =\frac{1}{2} \|{x}-{x_b} \|^2_{          { B^{-1}}} +\frac{1}{2}\|          y-H{(x)} \|^2_{           { R^{-1}}}  \  /            x \in D \}
\label{first}
\end{equation}
where $D$  is a grid domain in $\mathbb R ^{3}$ . In equation (\ref{first}) the vector $ {x}=\big[T,S,\eta,u,v \big]^\top $ is {an ocean} state vector composed by the temperature $T$, the salinity $S$,  sea level $\eta$ and  horizontal velocity field $(u,v)$. The vector  ${x_b}$ is the background state vector,  achieved by numerical solution of an ocean forecasting model and { is an approximation  of  the "true"  state vector $ {x_t}$}. The difference between background ${x_b}$ and any state vector ${x}$ is denoted by $          {\delta x}$:
\begin{equation}
       {   {x_b}=           {x}+            {\delta x }}
\end{equation}
The vectors ${x}$ and ${x_b}$  are  defined on the same space called {\em physical space}. The  vector ${y}$ in (\ref{first}) is the  observational vector defined on a different space called  {\em observational space} and the function $H$ is a non linear  operator that converts values defined in the physical space to values defined in the observational space.
{An ocean} state vector ${x}$ is related to observations  $ y$ by means the following relation:
\begin{equation}
{y}=H({x})+            {\delta y}
\end{equation}
where  $         {\delta y}$ is an effective measurement error. {In (\ref{first})  the matrix $\mathbf{R}= \langle  {\delta y_t} \ {\delta y_t}^\top  \rangle $, with ${\delta y_t}={y_t}-H({x_t})$, is the observational error matrix covariance and it is
assumed generally  to be diagonal, i.e. observational errors are seen as statistically independent. The
 $\mathbf{B}= \langle  {\delta x_t}\  {\delta x_t}^\top \rangle $, with ${\delta x_t}={x_b -x_t}$, is the background error matrix covariances and 
 is never assumed to be diagonal in its representation.}\\ 
Problem (\ref{first}) is solved by minimizing the following   explicit  form of cost function $J(         x)$: 
\begin{equation} 
J( x)=\frac{1}{2}({x}-{x_b})^{\top}\mathbf{B}^{-1}({x}-{x_b})+\frac{1}{2}(          y - H({x}))^{\top}{ \mathbf R}^{-1}(           y-H({x})).
\label{second}
\end{equation}

\noindent It is often numerically
convenient to exploit the weak non linearity of
$H$  by approximating  $H({x})$, for  small increments ${\delta x}$
,  with  a  linear approximation  around the  background vector ${x_b}$:

\begin{equation}
H(         x) \approx H({        x_b})+           \mathbf H           {\delta x}.
\label{fourth}
\end{equation}
where the linear operator  $\mathbf{H}$  is the $H$'s Jacobian  evaluated at $ {x}={x_b}$.  
\noindent The cost function $J$, using (\ref{fourth}), is approximated  by the following quadratic function:
\begin{equation}
J({  {\delta           x}})=\frac{1}{2}  { { {\delta             x}}}^{\top}\mathbf{B}^{-1} { {{\delta            x}}} +\frac{1}{2}(            d-         \mathbf  H{  {\delta           x}} )^{\top}\mathbf{R}^{-1}(            d -          \mathbf H {  {\delta            x}})
\label{functional2}
\end{equation}     
\noindent defined on increment space. In (\ref{functional2}) the vector $           d=          y -H(      x_b) $ is the misfit.\\
The minimum of the cost function $J(         { \delta x} )$ on the increment space may be justified by posing  $\nabla J({\delta x})=0 $. Then we 
obtain, as also shown in \citet{Haben},  the following  preconditioned system:


\begin{equation}
\bigg( \mathbf{I+BH^TR^{-1}H} \bigg)  {\delta x} =\mathbf{BH^TR^{-1}}         {d} 
\label{eight}    
\end{equation}

\noindent To solve the linear equation system (\ref{eight})  iterative methods able to converge toward a
practical solution are needed. Generally,{  the OceanVar model uses the Conjugate Gradient  Method \citep{Byrd}. } \\
The iterative minimizer schema is based essentially on matrix-vector operation  of some 
 vector $ {v}=\big(      \mathbf     {H^TR^{-1}H}\big)   {\delta           x}$ with the covariance matrix $\mathbf  B$. 
This computational kernel is required at each iteration and its huge computational complexity
is a bottleneck in practical data assimilation. This problem can be overcome by decomposing 
the covariance matrix $\mathbf B$ in the following form:
\begin{equation}
\mathbf{B} = \mathbf{V} \mathbf{V}^T
\end{equation}


\noindent However due to its still large size, the  matrix  $\mathbf{V}$ is {split}  at each minimization iteration as a sequence of
linear operators  \citep{Weaver}. {More precisely, in  OceanVar    the 
matrix $      \mathbf     V$ is decomposed as:}
\begin{equation}
\mathbf{V}=  \mathbf{V}_D \mathbf{V}_{uv} \mathbf{V}_{\eta}  \mathbf{V}_H \mathbf{V}_V
\label{composed}
\end{equation}
where the linear operator $\mathbf{V}_V$ transforms coefficients
which multiply vertical EOFs into vertical profiles of temperature 
and salinity defined at the model vertical levels, $    \mathbf       {V}_H$ and   $    \mathbf       {V}_\eta$ apply  respectively  the  gaussian filtering  to the  fields of temperature
and salinity,  and  sea surface.  $ \mathbf{V}_{uv}$  calculates velocity from sea surface height,
temperature and salinity, and $\mathbf V_D$ applies a divergence
damping filter on the velocity field. A more detailed   formulation of each linear operator is described in \citet{Dobricic}. In this paper we focus on 
operator $\mathbf {V}_H$.

\subsection{OceanVar Horizontal Covariances}
\noindent       The OceanVar  horizontal error covariances matrix   $  \mathbf   V_{H} $  is assumed  to be a  Gaussian matrix  \citep{Dobricic}.   In oceanographic  models, isotropic and Gaussian spatial correlations can
be relatively efficiently approximated  by an iterative application of a gaussian   RF  \citep{Lorenc,Hayden}  that requires only a few steps. Moreover, its application on a horizontal grid can be {split}  into two independent directions \citep{Purser}. We highlight that the ocean recursive filter scheme is more complicated {than}  the atmospheric case due to the presence of coastlines. In this framework, the horizontal  error covariances $\mathbf  V_H $ is factored as:
\begin{equation}
\mathbf{V}_{H}= \mathbf{G_y  G_x }
\end{equation}
 where $     \mathbf      {G_x}$ and $   \mathbf         {G_y}$ represent the gaussian operators in directions x and y.\\
In the next section we present a optimal revised RF to compute an approximation of the image of the temperature and salinity fields  by means of matrices  $\mathbf{G_x}$   and $ \mathbf{G_y}$. 

\section{Recursive Filters  for a  3D-VAR  Assimilation  Scheme}
\noindent Because of the separability of the  two-dimensional (2D) Gaussian function (that is $e^{-(x^2+y^2)}=e^{-x^2} e^{-y^2}$), a 2D-RF can be obtained applying a 1D-RF on each row and column of the discrete domain. Then we will consider the  properties  of the RFs only on one dimensional \citep[e.g.,][]{Oppe}.\\ 
\noindent The application of a one-dimensional n-th order RF on a grid of $m$ points  requires  two   main steps:


\begin{eqnarray}
 p^k_i =\beta s^{k-1}_i +\sum_{j=1}^n\alpha_j   p^{k}_{i-j}   \qquad  i=n+1,m   :+1,\  k=1,..,K  \label{one1}  \\
s^{k}_i          =\beta    p^{k}_i           +\sum_{j=1}^n\alpha_j s^{k}_{i+j} \ \ \  \qquad  i=m,n+1  :-1, \  k=1,..,K \label{two2}
\end{eqnarray}

\noindent where $K$  is the total filter iterations number, $s^0$ is the input distribution, $p^k $ 
is the  $k$-th  output  vector of the forward procedure (\ref{one1})  and  $s^k$ is the $k$-th   output   vector of the backward procedure (\ref{two2}),  corresponding to the input distribution for the  $(k+1)$-th forward procedure.  At last,  
$\alpha_j,\  j=1,...,n$  and $\beta $ are the filter  smoothing  coefficients.\\
\subsection{ The  1st-RF and  3rd-RF in OceanVar}
\noindent In the previous  OceanVar scheme, it was implemented a 1st-RF  algorithm  of \citet{Hayden,Purser}
along  the  x and   y directions.  The revised  scheme  uses a 3rd-RF based on the works of \citet{Young, Vlietnew}. {OceanVar model computes for each grid point  the filter coefficients that  depend  on the correlation radius $R(x,y)$ and 
the grid distances $\Delta x$ and $\Delta y$. Then OceanVar allows to use  different   length-scale for the gaussian covariance functions. \\}
\noindent The  one dimensional  1st-RF version   is composed  by the following rules: 

\begin{eqnarray}
 p^{k}_1=\beta_1 s^{k-1}_1,\ \ \ \  \qquad  \qquad \qquad \qquad \qquad \qquad k=1,....,K     \\
 p^k_i          =\beta_i s^{k-1}_i           +\alpha_i  p^{k}_{i-1}  \ \ \           \qquad  i=2,m   : +1,  \quad k=1,....,K    \label{proc1-1}    \\[3mm]
s^{k}_m=\beta_m p^k_m, \ \ \ \ \ \qquad  \qquad \qquad \qquad \qquad \qquad k=1,....,K     \\
s^{k}_i          =\beta_i p ^{k}_i           +\alpha_i s^{k}_{i+1}       \ 	\   \qquad  i=m-1, 2  :-1,    \ \  k=1,...,K  \   \label{proc1-2}   
\label{proc1}
\end{eqnarray}

\noindent where the  parameter $\alpha_i, \beta_i \in (0,1)$ and  $\beta_i=(1-\alpha_i)$ are the filter coefficients at the $i-th$ grid point. In order to obtain the filter smoothing coefficients $\alpha_i$ and $\beta_i$, the crucial  relationship in \citet{Hayden}  is considered:
\begin{equation}
R_i^2=2K\frac{\alpha_i}{(1-\alpha_i)^2} \Delta x_i^2.
\label{crucial}
\end{equation}
where $R_i$ and $\Delta x _i$ are  respectively  the correlation radius and  the grid distance  at the $i-th$ grid point.  \noindent  By  means of the equation (\ref{crucial}),  it follows that:
\begin{equation}
 (1-\alpha_i)^2 =2K   \big(\Delta x_i^2/ R_i^2   \big)  -  2K \big(\Delta x_i^2/ R_i^2  \big)    (1-\alpha_i). 
\label{evi}
\end{equation}
Calculating the roots  $\beta_i=(1-\alpha_i)$ from the equation (\ref{evi}),
  we obtain that   $\alpha_i$ and $\beta_i$ are: 
\begin{eqnarray}
\alpha_i=1+ \frac{K  \Delta x_i^2 }{ R_i^2 }  -\sqrt{  \frac{K \Delta x_i^2} { R_i^2}  \bigg(\frac{K  \Delta x_i^2} { R_i^2 } +2\bigg)}\\
\beta_i  = - \frac{K \Delta x_i^2}{ R_i^2 }  +\sqrt{  \frac{K  \Delta x_i^2} { R_i^2}  \bigg(\frac{K  \Delta x_i^2} { R_i^2 } +2\bigg)} \quad  
\end{eqnarray}


\noindent In the following  we are considering a 3rd-RF that  approximates  quite  successfully in   just one iteration   the gaussian convolution  considered in  \citet{Mirouzea} and  that allows different length scales for each grid point of computational domain .  It is based   on the works of \citet{Young} and \citet{Vliet} and  is  widely used for the blurring in  digital image and never used in geophysical data assimilation.\\
\noindent In the next theorem  we are  presenting a simple and accurate adapted version of this 3rd-RF   for  OceanVar scheme in   a one-dimensional
grid.
\begin{theorem}
For each $i-th$ grid point of a finite one-dimensional grid
  with  correlation radius $R_i$ and   grid spacing $\Delta x_i $,   a normalized 3rd-RF is given by:

\begin{eqnarray}
p_i=\beta_i s^0_i+\alpha_{i,1}  p_{i-1}+\alpha_{i,2} p_{i-2} +\alpha_{i,3} p_{i-3} \quad  i=4, m \ : +1 \label{fort} \quad  \ \  \\[2mm]
s_{i}=\beta_i p_i+\alpha_{i,1}s_{i+1}+\alpha_{i,2}s_{i+2} +\alpha_{i,3}s_{i+3}. \quad  i=m-3,1 \ : -1. \label{bact} 
\end{eqnarray} 
\noindent where the function  $s^0_i$ is the input distribution,  $\alpha_{i,1}=a_{i,1}/a_{i,0} $,  $\alpha_{i,2}=a_{i,2}/a_{i,0} $,  $\alpha_{i,3}=a_{i,3}/a_{i,0} $ and $\beta_i= \sqrt[4]{{2 \pi \big({R_i}/{\Delta x_i}\big)^2}} \big(1 - \big(\alpha_{i,1} +\alpha_{i,2}+\alpha_{i,3})\big)$ are the filtering  coefficients  and:
\small
\[
\begin{array}{l}
a_{i,0} = 3.738128   + 5.788982  \bigg( \dfrac{R_i} {\Delta x_i}\bigg)  + 3.382473  \bigg( \dfrac{R_i} {\Delta x_i}\bigg)^2 + 1.000000  \bigg( \dfrac{R_i} {\Delta x_i}\bigg)^3 \\[2mm]
a_{i,1}=  5.788982  \bigg( \dfrac{R_i} {\Delta x_i}\bigg)+6.764946  \bigg( \dfrac{R_i} {\Delta x_i}\bigg)^2+ 3.000000 \bigg( \dfrac{R_i} {\Delta x_i}\bigg) ^3 \\[2mm]  
 a_{i,2}= -3.382473 \bigg( \dfrac{R_i} {\Delta x_i}\bigg)^2-3.000000 \bigg( \dfrac{R_i} {\Delta x_i}\bigg)^3 \\[2mm]

a_{i,3}=  1.000000  \bigg( \dfrac{R_i} {\Delta x_i}\bigg)^3 \\
\end{array}
\]
\normalsize
\end{theorem}

{\footnotesize
\begin{proof}
 The procedure to obtain the 3rd-RF coefficients is given in Appendix A.   
$ $ \hfill $\blacksquare$
\end{proof}
}

\normalsize
\noindent By the Theorem 3.1  the   filtering strategy is the following:
\begin{itemize}
\item the  input data $s^0_i$ are first filtered in the forward direction as suggested by the difference equation in   (\ref{fort}) to get $p_i$.  
\item The output of this result $p_i$  is then filtered in the backward direction according to the difference equation corresponding to backward equation in (\ref{bact}) in order to get $s_i$.   
\end{itemize}

\noindent  To complete the statements  in (\ref{fort}) and (\ref{bact}) we  fix  the following heuristic initial conditions for
 the forward  and backward procedures:

\begin{equation}
\begin{array}{l}
forward\ conditions\\[2mm]
p_1=\beta_1 s^0_1    \\     
p_2=\beta_2 s^0_2+\alpha_{2,1}  p_{1}  \\
p_3=\beta_3 s^0_3+\alpha_{3,1}  p_{2}+\alpha_{3,2} p_{1}  \\[2mm]
\end{array}
\quad 
\begin{array}{l}
backward \ conditions\\[2mm]
s_m=\beta_m p_m    \\     
s_{m-1}=\beta_{m-1} p_{m-1}+\alpha_{m-1,1}  s_{m}  \\
s_{m-2}=\beta_{m-2} p_{m-2}+\alpha_{m-2,1}  s_{m-1}+\alpha_{m-2,2} s_{m}  \\[2mm]
\end{array}
\end{equation}

\noindent  Now we  give some remarks  on the accuracy  and the computational costs  of a  RF.

\begin{remark}[Remark 1]
{\em
For an arbitrary input  distribution $s^0_i$,  
a measure for   accuracy   of  a   RF is given by the following inequality:  
\begin{equation} 
\|            \epsilon_{s_i}  \|_2 \leq \|           \epsilon_{h_i} \|_2    \|          s^0_i  \|_2 
\label{ineq}
\end{equation}
where:
\begin{itemize}
\item  {$\|            \epsilon_{s_i} \|_2$ is  the euclidean norm of  the difference 
between the discrete  convolution $s^*_i=g_i  \otimes s^0_i  $(with $g_i$ normalized  gaussian  function) and the function   $s_i$, obtained by the   RF applied to $s^0_i$.}
\item   { $\|            \epsilon_{h_i} \|_2$   is    the euclidean norm of  the difference 
between the  gaussian  $g_i$ and   the function  $h_i$ (called \emph{impulse response}),   obtained by a RF applied to  Dirac  function;} 
\item   { $\| s^0_i\|_2$ is the euclidean norm  of the input function $s^0_i$ .  }
\end{itemize}
}

\end{remark}

{\footnotesize
\begin{proof}
The proof is shown in Appendix A
$ $ \hfill $\blacksquare$
\end{proof}
}

\noindent   By the  Remark  1  we observe that  the approximation error   $\|            \epsilon_{s_i}\|_2$ of $s^*_i$,  
  is  smaller than arbitrary  $\epsilon >0$ if and only if     
 $\|            \epsilon_{h_i}\| _2 < \epsilon / \|          s^0 \|_2 $. Then   to get a good approximation  of the gaussian convolution $s^*_i$, it is necessary that  the RF  determines a good approximation   $ h_i$ of  gaussian function  $ g_i$.

\noindent We highlight some practical considerations about the convergence of RF, applying it to a Dirac function . We choose  a  one-dimensional grid   of $m=300$ points , a constant  correlation radius $R= 12000 \ \emph m $ and a constant grid space $\delta x =6000 \ \emph m$.   In  Figure 1 we have the impulse response  $h_i  $, obtained by the  1st-RF and 3rd-RF, to reconstruct a Gaussian  function with  non-dimensional length scale  $\sigma=R/\delta_x$. In the left panel of Figure  1  we observe the slow convergence of   the 1st- RF   that needs an  high number of iterations to reach a good approximation of  gaussian function $g_i$. 
In the right panel of Figure  1, we  show the same  with the 3rd-RF. For this case, it is needed just one iteration  to obtain an  accurate  approximation $h_i$ of the gaussian function $g_i$. 
\begin{figure}[h!]
	\includegraphics[angle=0,width=0.52\linewidth]{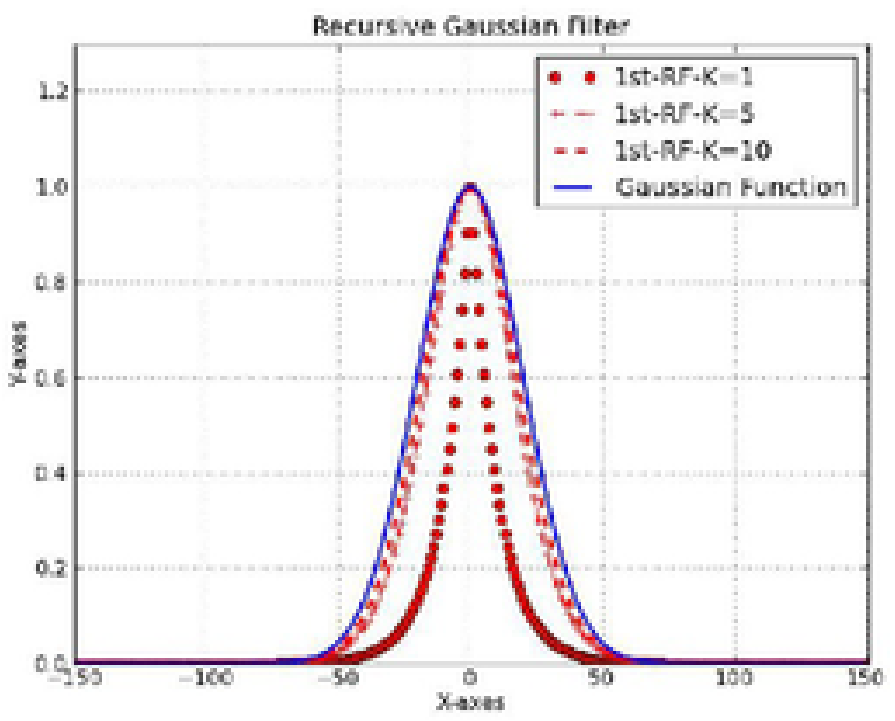}
       \includegraphics[angle=0,width=0.52\linewidth]{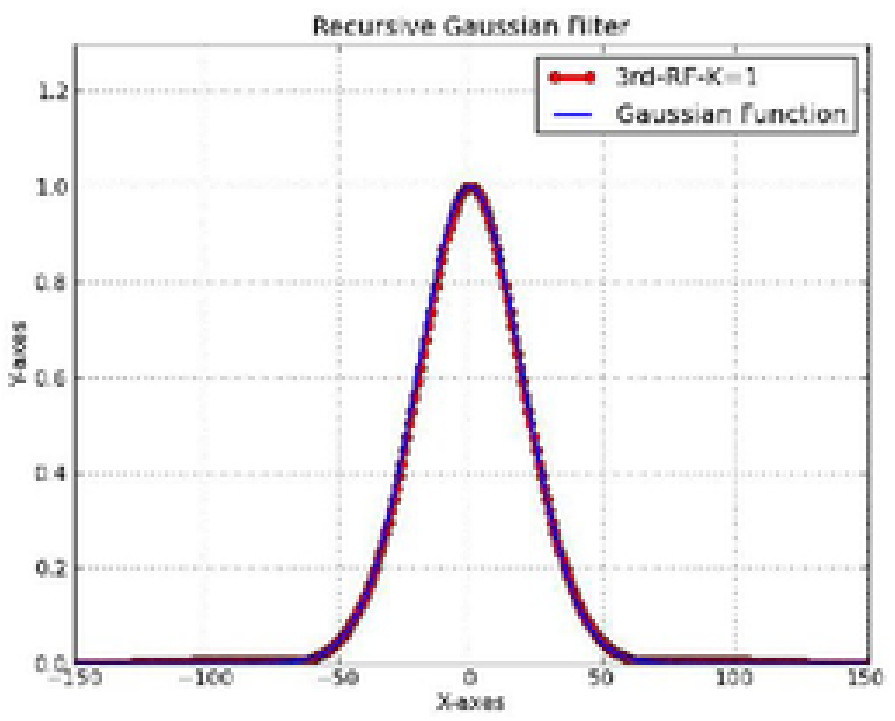}
\caption{\small{Left  -  \emph{The Impulse response  $h_i  $ (red) resulting from the application of equations    (\ref{proc1-1}) and (\ref{proc1-2})   for  $K=1,5,10$  and the true gaussian function $g_i$ (blue) with   non-dimensional length scale $\sigma=R/\delta _x$ and mean $\mu=0$ }\\
  Right - \emph{The Impulse response  $h_i $  (red) resulting from the application of equations    (\ref{for}) and (\ref{bac})   for  $K=1$  and the true gaussian function $g_i$ (blue) with   non-dimensional length scale $\sigma=R/\delta _x$ and mean $\mu=0$ } }}
\end{figure}

\noindent About the computational cost  of RF  we  consider the following Remark 2.
\begin{remark}[Remark 2]
{\em
The  computational time  of a  $n-$th order  accuracy RF  is  given by the following formula:

\begin{equation}
T(n,K,m)\approx  2\ (2n+1) \ m \ K   \ t_{calc}.
\label{calculus}
\end{equation}
where  $t_{calc} $ is the time for  a  floating point operation.
}

\end{remark}
 
\noindent It immediately follows that the time complexity of the 1st-RF is:
 \begin{equation}
T(1,K,m) \approx 6  m  K   \ t_{calc}
\label{storto1}
\end{equation}
\vspace{1mm}

\noindent Unfortunately, the  1st-RF needs a high step number $K$ to determine a good approximation of the convolution function. We underline that this is a huge obstacle in the OceanVar framework  since the RF iterations are too much {expensive} from a computational point of view for  real applications. Although the 3rd-RF complexity is
\begin{equation}
{T(3,K,m) \approx 14 \ m K\    \ t_{calc}}
\label{storto2}
\end{equation}

\noindent we have good approximation of the convolution kernel {within} just one iteration. Note that by comparing the time complexities given by (\ref{storto1}) and (\ref{storto2}) it follows that the theoretical computational time of the 1st-RF is less than that of the 3rd-RF only if 1 or 2 iterations are used, which however provides a very inaccurate approximation of the Gaussian function.  \\In the next section, we prove it through numerical experiments, confirming that the new 3rd-RF  in OceanVar can  improve the entire performance of the data assimilation  software.


\section{ Experimental  Results  of OceanVar  using 1st-RF  and 3rd-RF in Mediterranean Sea Implementation}
\noindent {In this section we test  the OceanVar set-up  in the Mediterranean Sea with  both  
 1st-RF and 3rd-RF from a numerical point of view. Either  1st-RF and 3rd-RF allow to use  different filtering scales as we will show in the next section  for global applications but in this configuration we fix a constant length scale.}   This  set-up follows the configuration  of \citet{Dobricic}, but the horizontal resolution is two times higher. In particular, the model has 72 horizontal  levels, and
the  horizontal resolution  is about 3.5 km in the latitudinal direction and between 3 km and 2.5 km in the longitudinal direction with an horizontal  grid of   $(1,742, 506)$ points.   The Mediterranean Sea  has a relatively large variability of the bathymetry, with both deep ocean basins, like the Ionian Sea, the Levantine and the Western Mediterranean, and extended  narrow shelves.\\     
\noindent  In  this configuration  we  apply  the  1st-RF  and the  3rd-RF  with the assimilation of only
 Argo floats profiles \citep{Poulain}.
We compare  accuracy, computational  time  and  the memory usage. The numerical tests are carried out on a parallel IBM cluster (with  {30 IBM P575 nodes, 960 cores Power6 4.7GHz,  Infiniband 4x DDR interconnection, Operating system  AIX v.5.3}).

\noindent  In the following  we  start to compare the accuracy  of  horizontal   covariances in temperature by means of the two different RFs.\\
In Figure 2,    we report analysis increment obtained using only one iteration of 3rd-RF    on the Mediterranean Sea, for the temperature   at $300m$, assuming
isotropic and Gaussian spatial correlations by means a constant correlation radius ($R = 15000m$).   The impression in Figure 2 is that the 3rd-RF reconstructs quite successfully  the horizontal  covariances in just one iteration. {The   long range correlation shows more  a diamond-like  shape but the inner part, the most intense  part of the field,  has the right shape and amplitude.} {Moreover we 	
underline that the 3rd-RF models the horizontal covariances  near the costs using the same number of the imaginary sea points used for the 1st-RF \citep{Dobricic} because the approximate gaussian functions built by both  the RFs have the same length scales}.\\
\noindent  Then  we compute  the horizontal  temperature   covariances at same depth  by means of the  1st-RF using K=1 iteration (1st-RF-K=1) in Figure 3a, 1st-RF,  K=5 iterations (1st-RF-K=5) in Figure 3b and 1st-RF, K=10 iterations (1st-RF-K=10) in Figure 3c.  We can  note that  only  with the 1st-RF K=10 (Figure 3c)  the horizontal  covariances became  isotropic and Gaussian. 
 Figure 4  shows the absolute differences between the horizontal temperature  covariances  computed by 1st-RF-K=1,5,10  and  the  3rd-RF. Observing the maximum values ​​of the differences in the three case, it is evident   that  the same accuracy of the 3rd-RF can be achieved only with a large number of iterations of the other one . Qualitatively the results  for the salinity are the same as for  temperature and hence we do not show them.
\begin{figure}[ht!]
			\centering
                    \includegraphics[scale=.42]{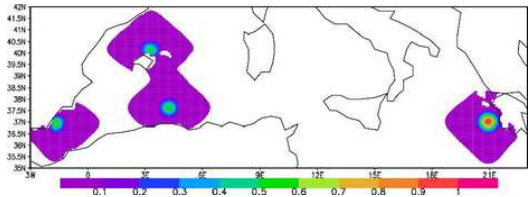} \\[3.5mm]  
         \caption{ \emph{{Analysis increments of temperature at 300 m
of depth for the Mediterranean Sea configuration  using 3rd-RF.}}}
\end{figure}

%
\begin{figure}[ht!]
			\centering   
                     \includegraphics[scale=.42]{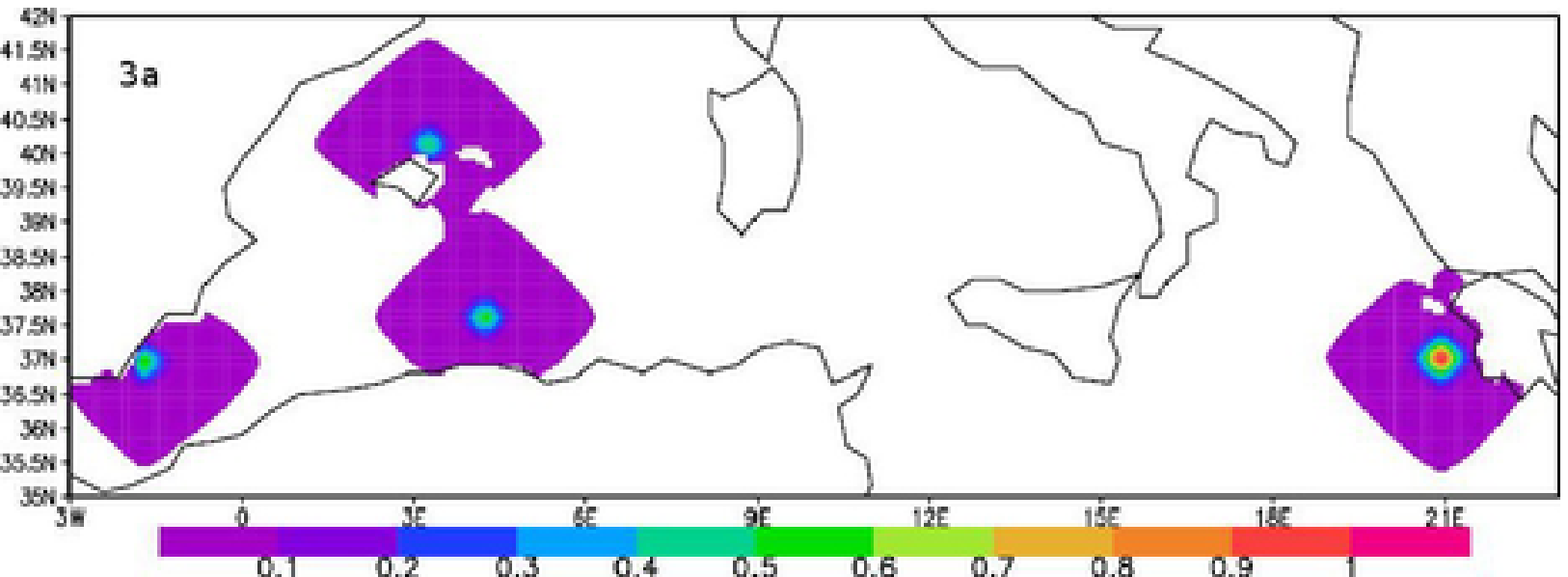} \\[3.5mm]
                    \includegraphics[scale=.42]{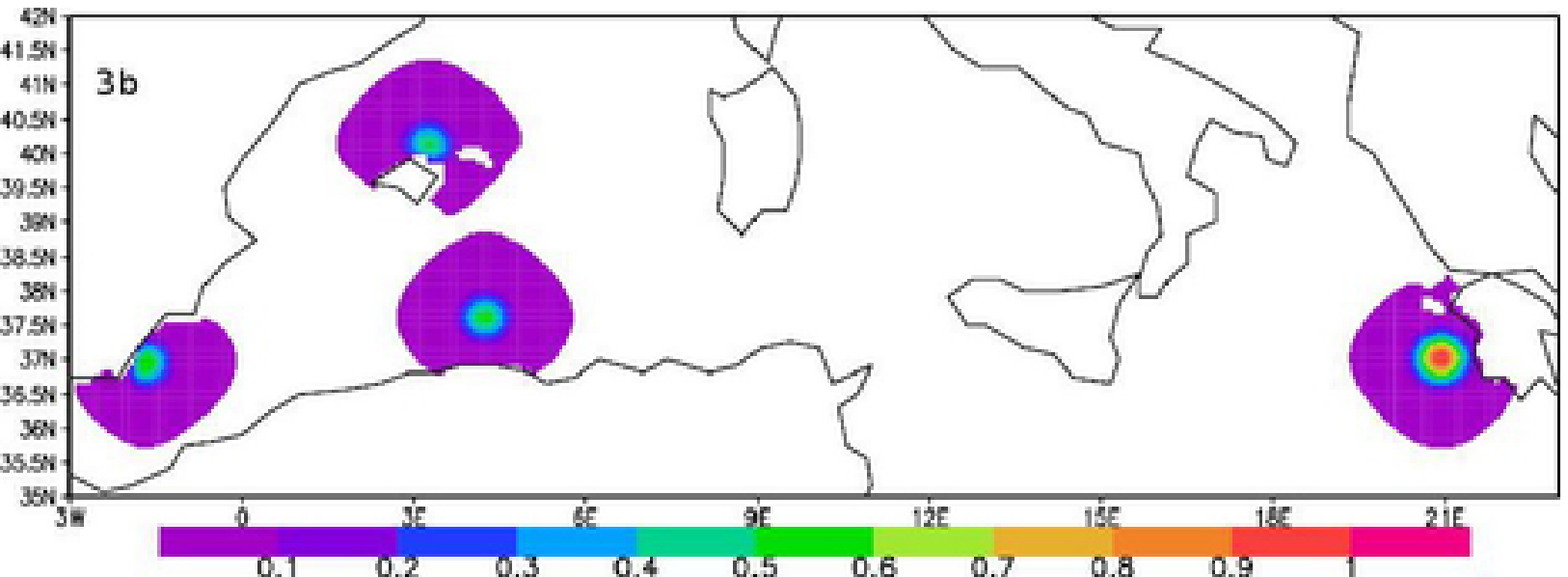}\\[3.5mm] 
                     \includegraphics[scale=.42]{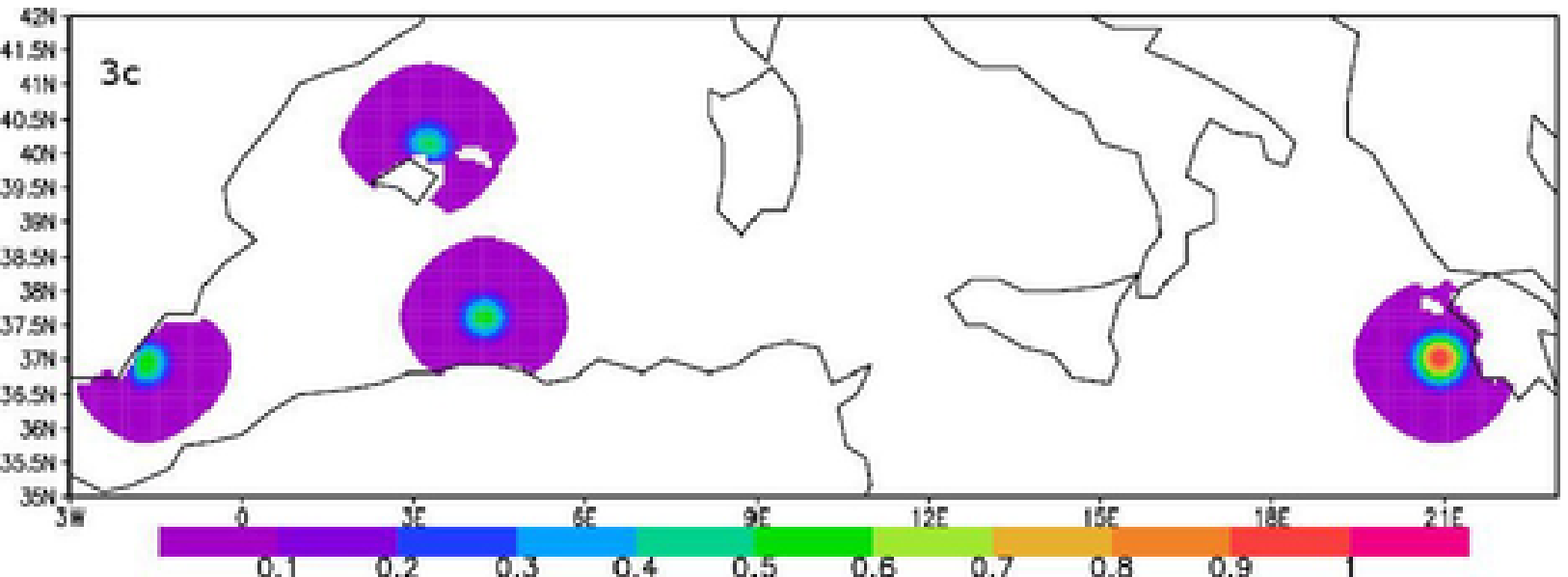} \\[3.5mm]  

			\caption{ \emph{{Analysis increments of temperature at 300 m
of depth for the Mediterranean Sea configuration  using 1st-RF-K=1 (3a),    1st-RF-K=5 (3b)
   and   1st-RF-K=10 (3c).}}}
\label{tpipe1}
\end{figure}

\begin{figure}[ht!]
			\centering
                     \includegraphics[scale=.41]{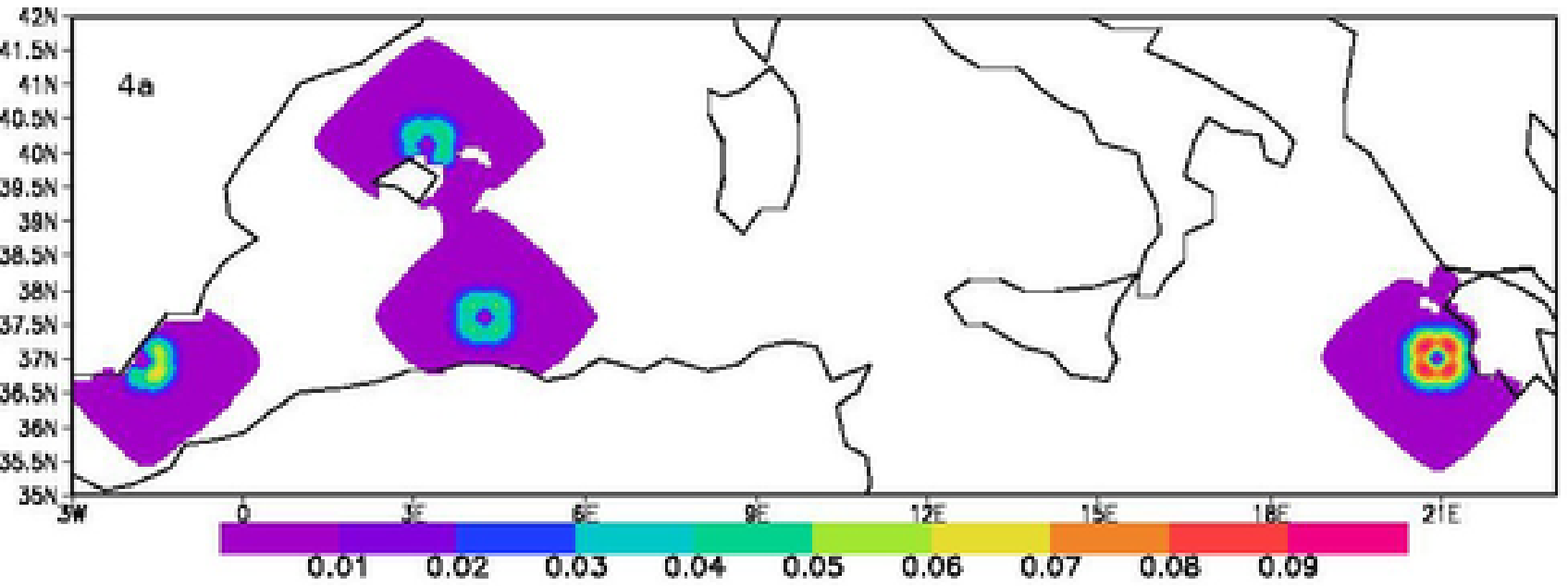} \\[3.5mm]
                    \includegraphics[scale=.41]{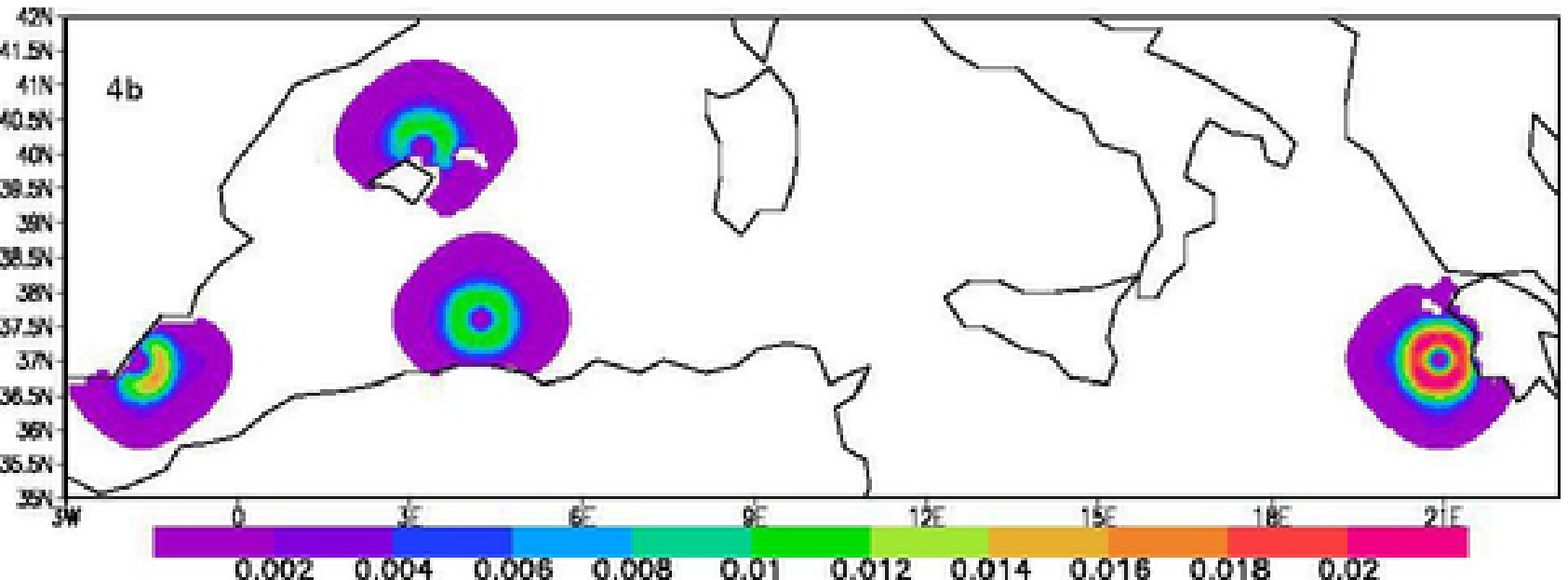} \\[3.5mm]
                    \includegraphics[scale=.41]{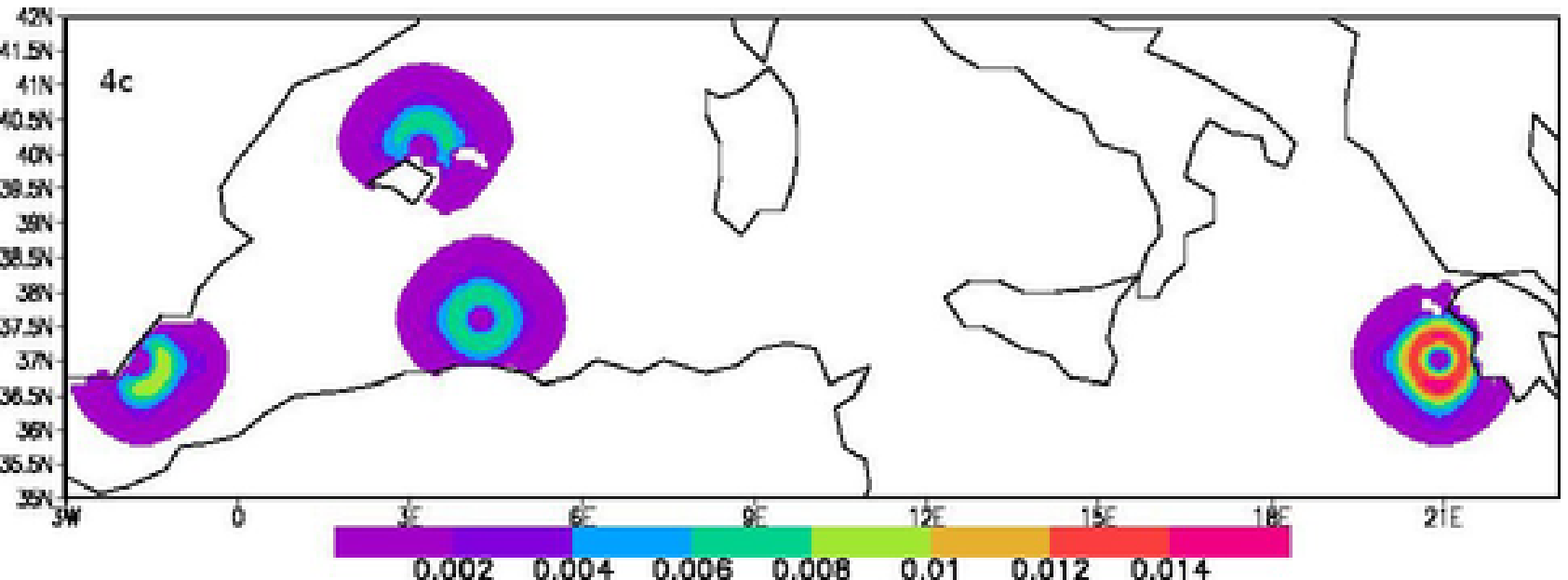}    
                	\caption{ \emph{{Absolute Differences for the temperature at 300 m of depth for the Mediterranean Sea configuration between 1st-RF-K=1 and  3rd-RF (4a), 1st-RF-K=5 and  3rd-RF (4b), 1st-RF-K=10 and  3rd-RF (4c). }}}
\label{tpipe2}
\end{figure}


\subsection{ Performance Results in Mediterranean Sea Implementation }
\noindent In  previous section  we show that  the  use of a more  accurate  RF  needs less iterations in order to obtain isotropic and Gaussian spatial
 horizontal  covariances.\\
In this subsection we  report the benefits and disadvantages of  the 3rd-RF  from a computational point of view.
In particular here  we give  some information on the execution time  and the  memory  usage in  the parallel version of the 1st-RF and  3rd-RF of Oceanvar in Mediterranean Sea  configuration.\\ 
  The parallel implementation of the RF in Oceanvar  uses  communication strategies  based on the pipeline method \citep[e.g.,][]{Nakano}, because   RF is a typical algorithm
with flow dependences,  where  each iteration has to be strictly executed in a pre-fixed order. OceanVar  in Mediterranean Sea  configuration implements the pipeline method  for the RF by using  a   {column-row-wise block distribution of processors} and  {blocking {\tt mpi-send} and {\tt mpi-receive} communication functions  to transfer  the boundary conditions. All the other operators in OceanVar (vertical EOFs, Vv, barotropic operator, Vn, velocities operator, Vuv, and divergence damping filter, Vd), along with the computations in observation space (misfits update, etc.) follow a Cartesian domain decomposition. The parallel implementation is also completely independent from that of the ocean model used for the forecasts.
We modify OceanVar  by  substituting the 1{st}-RF with the 3{rd}-RF,  without modifying  {the pipeline method}. \\
  As  3rd-RF   needs only one  iteration, it requires  a smaller number of mpi-communications. This results in  a sensible reduction in the  latency in the communications. However  3rd-RF has  to transfer  a larger amount of data in each iteration.  
For example in the case of 1st-RF, in the forward filtering,  the left processor   sends only   data on  its last column  to the right processor  while in the case of   3rd-RF,  
the left processor sends   data on  its last  three columns  to  the right processor.\\
 Figure 5  shows  the  execution time in seconds   and the memory usage in Gbyte  of the   OceanVar   with  1st-RF-K=1, 1st-RF-K=5, 1st-RF-K=10 and 3rd-RF on the IBM cluster using 64 processors.  In particular  we   report the  OceanVar   and the RF wall clock  times. To estimate the RF times  in the software  we   synchronize all the processes in our  mpi-communicator  by  using  MPI-Barrier, before and after the implementation of RF. \\
\noindent Figure 5  shows that   Oceanvar with  the  3rd-RF gives the best results in terms of  the  execution time if we consider  the accuracy on the horizontal   covariances discussed in previous section. The 3rd-RF compared to the 1st-RF-K=5 and  the 1st-RF-K=10 reduces the wall clock time of the software  respectively of about $27\%$  and  $48\%$. 
The execution time is reduced because the new RF decreases the sequential time of each process in the  mpi-communicator and needs  only one  {\tt mpi-send} and {\tt mpi-receive} communication for each process.
 The  3rd-RF algorithm uses twice as much memory   than 1st-RF,  because it has twice the number of the filter coefficients.  However, this did not represent a problem on our cluster  because  the maximum memory allowed was about 250 Gb. \\
\begin{figure}[h!]
			\centering
                       \includegraphics[scale=.7]{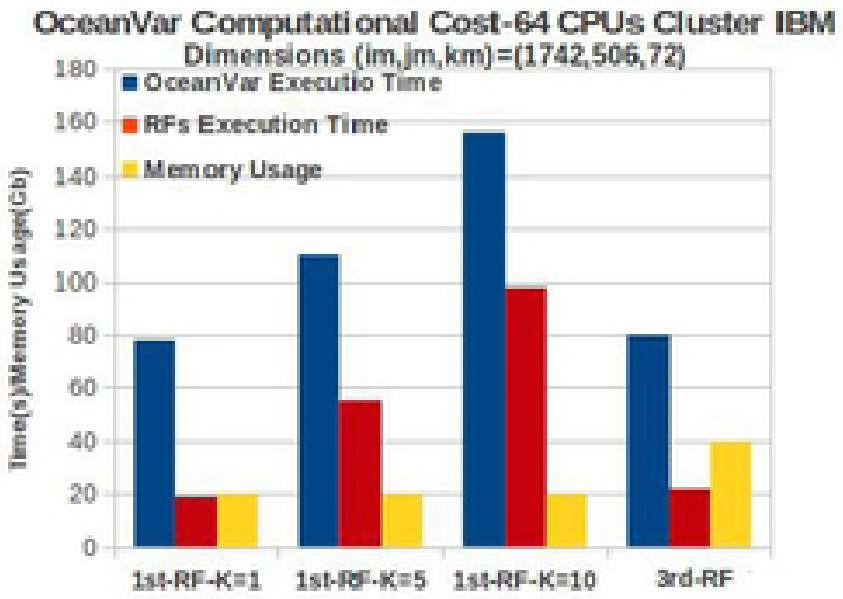}\\[2mm]
\caption{ \emph{{Performance results (execution time and memory usage) for the
recursive filter tests within the Mediterranea Sea implementation of OceanVar.
 }\\
 }}
\label{pipe}
\end{figure}  

\section{Experimental Results of OceanVar   using 1st-RF  and 3rd-RF  in a Global Implementation}
\noindent In this section we describe experimental results of the 3rd-RF in a Global Ocean implementation
of OceanVar that follows \citet{mwr}. The model resolution is about ${1}/{4}$ degree and the horizontal grid is  tripolar,
as described by \citet{tripolar}. This configuration of the model is used at CMCC for global ocean physical reanalysis  
applications \citep{ferry}.
The model has 50 vertical depth levels. The three-dimensional model grid consists of 73614100 grid-points. The comparison between the 1st-RF and
the 3rd-RF is here carried out for a realistic case study, where all in-situ observations of temperature and salinity from Expendable bathythermographs (XBTs), Conductivity, Temperature, Depth (CTDs) Sensors, Argo floats and Tropical mooring arrays are assimilated.
The observational profiles are collected, quality-checked and distributed by Coriolis \citep{cabanes}. The global application of the recursive filter
accounts for spatially varying  and season-dependent correlation length-scales (CLSs), unlike the Mediterranean Sea implementation.
Correlation length-scale were calculated by applying the approximation given by \citet{bpb}
to a dataset of monthly anomalies with respect to the monthly climatology, with inter-annual trends removed.
The two panels of Figure \ref{fig_glob1} shows an example of zonal and meridional temperature correlation length-scales
at 100 m of depth, respectively, for the winter season in the Western Pacific.
 \begin{figure}[t]
\centering
\includegraphics[width= 1.\textwidth,clip,angle=0]{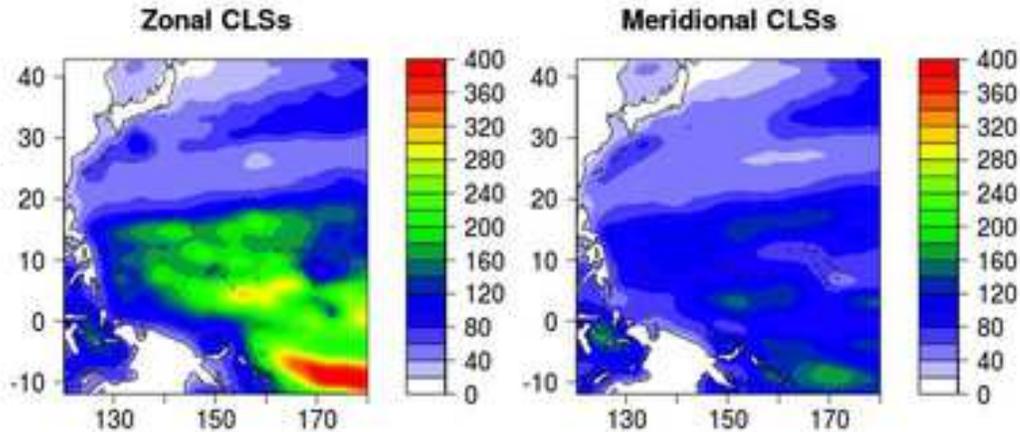}
\caption{\emph
{Zonal (left) and meridional (right)
correlation length-scales for temperature
at 100 m of depth for the Western Pacific Area.
}}
\label{fig_glob1}
\end{figure} Typically, areas characterized by strong variability
(e.g.~Kuroshio Extension) present shorter correlation length-scales of the order of less than 100 Km that lead to very narrow corrections,
while in the Tropics the length-scales are longer and may reach up to 350 Km, thus broadening
the 3DVAR corrections. Furthermore, at the Tropics, it is acknowledged that
zonal correlations are longer  than the meridional ones \citep[e.g.][]{derbros}. The analysis increments from a 3DVAR applications
that uses the 1st-RF with 1, 5 and 10 iterations and the 3rd-RF are shown in Figure \ref{fig_glob2},
with a zoom in the same area of Western Pacific Area as in Figure \ref{fig_glob1}, for the temperature
at 100 m of depth. \begin{figure}[h!]
\includegraphics[width=.92\textwidth,clip,angle=0]{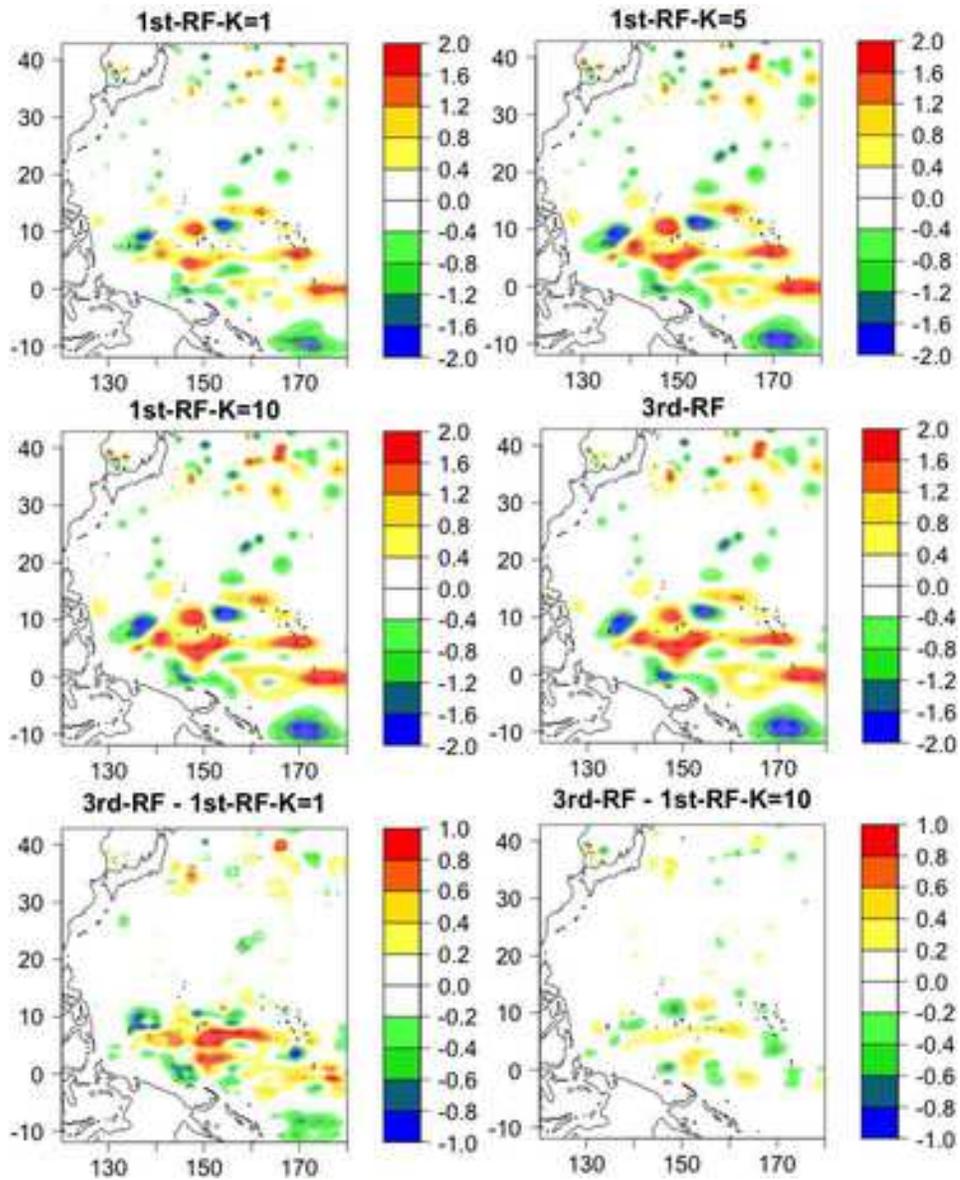}
\caption
{Analysis increments of temperature at 100 m
of depth for the Western Pacific for different
configurations of the recursive filter (first two rows of panels).
Differences of 100 m temperature analysis increments
between 3rd-RF and 1st-RF-K=1 and between 3rd-RF and 1st-RF-K=10 (bottom panels).
}
\label{fig_glob2}
\end{figure} The figure also displays the differences between the 3rd-RF and the 1st-RF with either 1
or 10 iterations. The patterns of the increments are closely similar, although increments for the case of 1st-RF-K=1
are generally sharper in the case of both short (e.g.~off Japan) or long (e.g.~off Indonesian region) CLSs. The panels of the differences reveal also
that the differences between 3rd-RF and the 1st-RF-K=10 are very small, suggesting once again that the same
accuracy of the 3rd-RF can be achieved only with a large number of iterations for the first order 
recursive filter.

\subsection{Performance Results in Global Implementation}
\noindent In this subsection we present the performance results of the 3rd-RF with respect to the 1st-RF for the
Global Ocean case study.  OceanVar is run on a
parallel IBM cluster, each of the 482 nodes with two eight-core Intel Xeon processors.
We use a cluster different  from the one  previously introduced for the Mediterranean Sea configuration, in an attempt
of presenting performance results also in different computing environments.
The global ocean implementation of OceanVar has a parallel implementation that differs
from the pipeline method previously presented: it exploits hybrid MPI-OpenMP parallelism, where OpenMP
acts over the vertical level loops while MPI over the horizontal grid.
This strategy has the advantage of limiting the MPI communication when exploiting the same
number of cores, with the shortcoming of being not very flexible (due to the limitations of
the cores-per-node number of threads for the OpenMP vertical parallelism).
We have used 5 nodes for our tests: 5 MPI processes, 16 threads for a total of 80 cores used.
Results are summarized in Figure \ref{fig_glob3}, when the number of iterations of the 3DVAR minimizer
is fixed to 30, as in realistic applications.\begin{figure}[t]
\centering
\includegraphics[width=.90\textwidth,clip,angle=0]{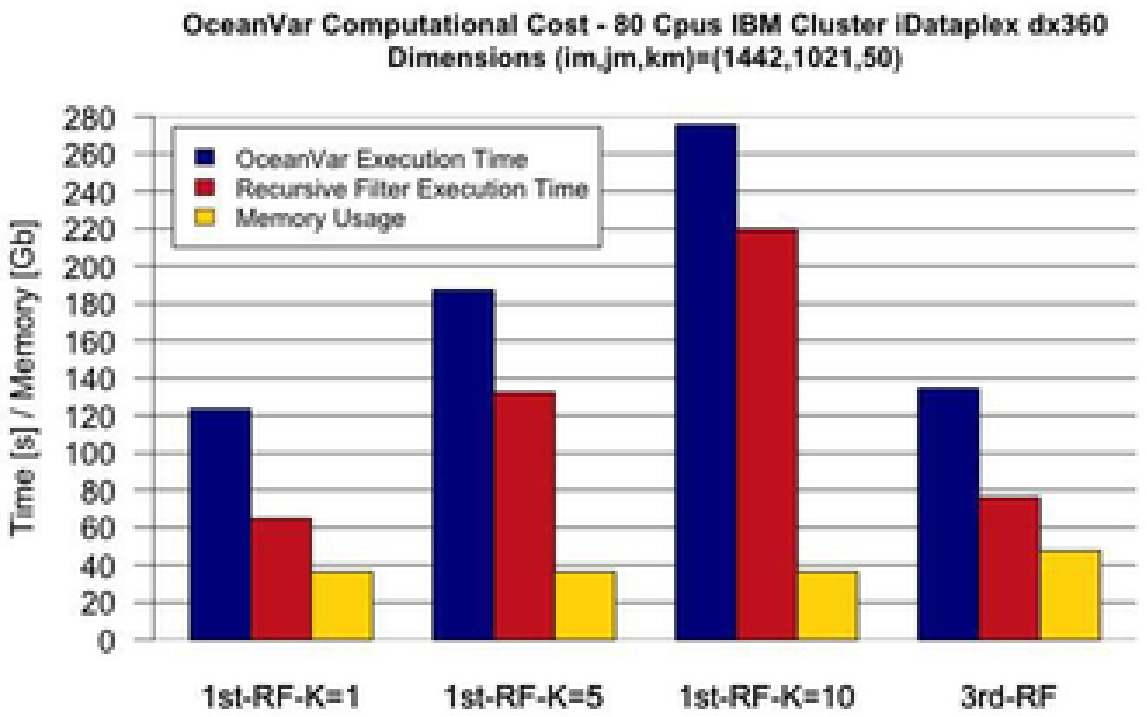}
\caption
{Performance results (execution time and memory usage)
for the recursive filter tests within
the Global Ocean implementation of OceanVar.
}
\label{fig_glob3}
\end{figure}  Generally, relative performances of the two recursive filters are comparable with those
of the Mediterranean Sea implementation. 3rd-RF reduces the total wall clock time
of OceanVar by  28 $\%$ and $51\%$ with respect to  1st-RF-K=5 and 1st-RF-K=10, respectively.
This reduction increases up to 42 $\%$ and 65 $\%$ if we consider only the recursive filter routines.
On the other hand, the total memory usage increases from 36.4 Gb (7.3 Gb per node) to
47.6 Gb (9.5 Gb per node), i.e.~by  30 $\%$. Thus, the third order recursive filter is able to significantly reduce the execution
time at the price of an affordable increase of memory usage.

\section{Conclusions}
\noindent In this paper we describe the development and implementation of a revised scheme to compute the horizontal    covariances of temperature and  salinity  in an  oceanographic variational  scheme.\\ 
The existing  recursive filter (RF) of the first order (1st-RF)  is substituted with a  recursive filter of  the third order (3rd-RF). {Numerical experiments in Mediterranean Sea and Global Ocean demonstrated that a 3rd-RF  can significantly reduce the total computational time of the data assimilation scheme  maintaining the same level of accuracy}.\\
In addition  we provide the full theoretical development of the new filter,  based on the study by \citet{Young} and \citet{Vliet},  who  formulated the filter in the context of signal processing. We  adapt it  for a 3D-VAR oceanographic  scheme with the detailed description of the process  to obtain the 3rd-RF  coefficients for different length scales.  .\\ 
{Our implementation  is different  to that by \citet{Purser}, since  we have applied a different mathematical methodology to calculate the filter coefficients}.  The new RF is faster than the previous one because it requires only one iteration to compute the horizontal covariances. Therefore, we may assume that it will be faster than other first order accurate methods that need several iterations like the one described by \citet{Mirouzea}. This hypothesis  is  tested  with the pipeline method of Mediterranean Sea implementation  and  hybrid MPI-OpenMPI parallelization strategy  of the Global Ocean configuration presented in the previous sections. By applying some other parallelization strategy the relative performance of the new filter may differ. However, we believe that other parallelization strategies  are overall much less efficient than  those presented.\\
 The  future  improvement of the  3rd-RF  scheme in OceanVar will be the implementation of different mathematical  boundary  conditions at the coasts  and  a formulation of  3rd-RF for  spatially inhomogeneus and anisotropic covariances. 

\section{Acknowledgment}\noindent The research leading to these results has received funding from the
Italian Ministry of Education, University and Research and the Italian Ministry of Environment, Land and Sea under the GEMINA  and Next Data projects.

\section{Appendix A}
\begin{remark}[Theorem 3.1]
For each $i-th$ grid point of a finite one-dimensional grid
  with  correlation radius $R_i$ and   grid spacing $\delta x_i $,   a normalized  3rd-RF is given by:

\begin{eqnarray}
p_i=\beta_i s^0_i+\alpha_{i,1}  p_{i-1}+\alpha_{i,2} p_{i-2} +\alpha_{i,3} p_{i-3} \quad  i=4, m \ : +1 \label{for} \quad  \ \ \\[2mm]
s_{i}=\beta_i p_i+\alpha_{i,1}s_{i+1}+\alpha_{i,2}s_{i+2} +\alpha_{i,3}s_{i+3}. \quad  i=m-3,1 \ : -1. \label{bac} 
\end{eqnarray} 
\noindent where the function  $s^0_i$ is the input distribution,  $\alpha_{i,1}=a_{i,1}/a_{i,0} $,  $\alpha_{i,2}=a_{i,2}/a_{i,0} $,  $\alpha_{i,3}=a_{i,3}/a_{i,0} $ and $\beta_i= \sqrt[4]{{2 \pi \big({R_i}/{\delta x_i}\big)^2}} \big(1 - \big(\alpha_{i,1} +\alpha_{i,2}+\alpha_{i,3})\big)$ are the filtering  coefficients  and:
\small
\[
\begin{array}{l}
a_{i,0} = 3.738128   + 5.788982  \bigg( \dfrac{R_i} {\delta x_i}\bigg)  + 3.382473  \bigg( \dfrac{R_i} {\delta x_i}\bigg)^2 + 1.000000  \bigg( \dfrac{R_i} {\delta x_i}\bigg)^3 \\[2mm]
a_{i,1}=  5.788982  \bigg( \dfrac{R_i} {\delta x_i}\bigg)+6.764946  \bigg( \dfrac{R_i} {\delta x_i}\bigg)^2+ 3.000000 \bigg( \dfrac{R_i} {\delta x_i}\bigg) ^3 \\[2mm]  
 a_{i,2}= -3.382473 \bigg( \dfrac{R_i} {\delta x_i}\bigg)^2-3.000000 \bigg( \dfrac{R_i} {\delta x_i}\bigg)^3 \\[2mm]

a_{i,3}=  1.000000  \bigg( \dfrac{R_i} {\delta x_i}\bigg)^3 \\
\end{array}
\]
\normalsize

\label{teorema1}
\end{remark}

{\footnotesize
\begin{proof}[Proof.]   
\noindent  In order to obtain an efficient  RF  such that approximates  the gaussian convolution (as considered in   \citet{Mirouzea}):
\begin{equation}
s(x)= g(x) \otimes s^0(x) =\int_{-\infty}^{+\infty} g(x-\tau)s^0(\tau) d\tau,
\label{convolution}    
\end{equation}
we apply  the  Fourier  transformation  to   the equation  in (\ref{convolution}), where  $g(x)=\exp(- \frac{x^2}{2 \sigma^2})$ is  the  normalized  gaussian function of mean $\mu=0$ and    non-dimensional length-scale  $\sigma=R_i/\delta x_i$  associated to i$-th$ grid point. Hence   we  obtain:
\begin{equation}
          {S(w)}  =           {G(w)}       {S^0(w)}    \quad w \in     \mathbb{R} 
\label{lapli}
 \end{equation}
where  $          {S(w)},           {G(w)}$   and   $          {S^0(w)}$ in (\ref{lapli})  are respectively  the Fourier   transformations of  the function  $  s(x),    g(x)$,  and   $s^0(x)$.  For the   Fourier  transformation  of $g(x)$,  we have the well-note result:
\begin {equation}
G(w) =\sqrt{2 \pi} \sigma e^{-\frac{(\sigma \omega)^2}{2}}.
\label{gau} 
\end{equation}  
 \noindent  Using  a rational  approximation of the gaussian function in \citet{Abramo}:
\begin{equation}
\frac{1}{\sqrt{2 \pi}} e^{-\frac{t^2}{2}}= \frac{1}{b_0 +b_2t^2 +b_4 t^4+b_6 t ^6} + \epsilon (t), \quad t \in \mathbb{R}
\end{equation}
where  $b_0=2.490895$, $b_2=1.466003$, $b_4=-0.024393$,  $b_6=0.178257$ and $\epsilon (t)< 2.7*10^{-3}$, then
we can approximate  the  $G(w)$ function as:     
\begin{equation}
G(w) \approx  H_\sigma(w) = \frac {{2 \pi } \sigma C^2}{b_0 +b_2( \sigma w)^2 +b_4  ( \sigma  w)^4+b_6 (  \sigma  w) ^6}
\label{riferimento}
\end{equation} 
where $C^2$ is a square  of  a normalization constant that we will choose later. The absolute error $| \epsilon(w)|$ between $G(w)$ and $ H_\sigma(w)$ is less than ${5.4 \pi } \sigma  \times 10^{-3}$.
Moreover we can rewrite the $H_\sigma (\omega)$ as a  function in the  complex field:
\begin{equation}
 H_\sigma(s) = \frac {{2 \pi} \sigma C^2} {b_0 -b_2( \sigma s)^2 +b_4  ( \sigma  s)^4-b_6 (  \sigma  s) ^6}, \ \ with \  s=i \omega \in \mathbb{C} 
\end{equation} 
 \noindent Through the  non linear Newton Formula    \citep[e.g.,][]{Quart}, we  determinate  two real solution  $s=\pm 1.16680481/ \sigma$ of  the polynomial  $b_0 -b_2( \sigma s)^2 +b_4  ( \sigma  s)^4-b_6 (  \sigma  s) ^6$. After we divide  it   for  $\big((\sigma s)^2-1.6680481^2 \big)$ and
  obtain a   polynomial quotient  that is   difference between  a polynomial of fourth degree  and  second one.  Hence the rational polynomial  $H_\sigma (s)$ can be  decomposed in the following way:
\begin{equation}
H_\sigma ( s)= H_{\sigma,b}(s)  H_{\sigma, f}( s)
\end{equation}   
where
\begin{eqnarray}
H_{\sigma,f}( s) =\dfrac{ \sqrt { (2 \pi  \sigma )} C  }{  (1.166805+\sigma s)   ( 3.203730+2.215668 \sigma s+ ( {\sigma s})^2) }   \label{indietro}  \\
H_{\sigma,b}( s) = \dfrac{ \sqrt { (2 \pi  \sigma )} C  }{   (1.166805-\sigma s)    (  3.203730-2.215668 \sigma s  +  ({\sigma s})^2)}   \label {avanti}    .
\end{eqnarray}
  
\noindent Using standard approximations as the  backward and forward  differences  for the  $\mathcal{Z}$  transformation   \citep[e.g.,][]{Oppe}   to switch  from continuous to discrete  
problem, we   replace $s=1-z^{-1}$  in (\ref{indietro}) and   $s=z-1$ in (\ref{avanti}) as in \citet{Young}. Hence we  obtain:

\begin{equation*}
H_{\sigma,f}( z^{-1}) =\dfrac{ \sqrt{ (2 \pi  \sigma)}  C }{  \big(1.166805+\sigma (1-z^{-1})\big)   \big( 3.203730+2.215668 \sigma (1-z^{-1})+ ( {\sigma (1-z^{-1})})^2\big) }
\end{equation*}
\begin{equation*}
H_{\sigma, b}(z) = \dfrac{ \sqrt{ (2 \pi  \sigma)}  C }{   \big(1.166805-\sigma (z-1)\big)  \big(  3.203730-2.215668 \sigma {(z-1})  +  ({\sigma (z-1)})^2\big)} 
\end{equation*}

\noindent  Both previous equations can be rewritten respectively  as standard  polynomials in  $z^{-1}$  and $z$:

\begin{eqnarray}
H_{\sigma,f}( z^{-1}) =\dfrac{ \sqrt{ (2 \pi  \sigma)} C}{  \big(a_{i,0}-a_{i,1} z^{-1} -a_{i,2} z^{-2}-a_{i,3}z^{-3} \big) }   \label{indietro2}   \ \   \\
H_{\sigma,b}(z) = \dfrac{{   \sqrt{ (2 \pi  \sigma)} C         }}{   \big( a_{i,0}-a_{i,1} z^{1} -a_{i,2} z^{2} -a_{i,3}z^{3} \big)   }   \qquad  \ \  \ \   \label {avanti2}           .
\end{eqnarray}

\noindent where  we can determinate the following coefficients for (\ref{indietro2}) and   (\ref{avanti2}): 

\begin{equation} 
\begin{array}{l}
a_{i,0} = \big ( 3.738128   + 5.788982 \sigma  + 3.382472 \sigma^2 + 1.000000 \sigma^3)\\[1mm]
a_{i,1}=     \big (5.788824 \sigma+6.764946 \sigma^2+ 2.999999\sigma ^3) \\  [1mm]
 a_{i,2}=   \big(-3.382472 \sigma^2- 2.999999\sigma ^3)  \\[1mm]
a_{i,3}=    \big(1.000000 \sigma^3)  
\end{array}
\label{angela}
\end{equation}
The normalized constant $C$  can be specified   by using the constrain that  the approximation gaussian  function $H_\sigma(\omega)$  must be $\sqrt{ 2 \pi } \sigma$   for  $\omega=0$  (see the equations (\ref{gau}) and (\ref{riferimento}))  hence $ H_{\sigma,f} ( z^{-1})=1 \wedge H_{\sigma,b} ( z)= \sqrt[4]{ (2 \pi  \sigma^2)}$ for $z^{-1}=1 \wedge z=1 $ . Follows that:

\begin{equation}
 C=\frac{\big( a_{i,0}-(a_{i,1}  +a_{i,2}  +a_{i,3}) \big)}{ \sqrt[4]{ (2 \pi)}}
\end{equation}

\noindent  Place  $ P(z)=H_{\sigma,f}(z)S^0(z)$   and  $S(z)=H_{\sigma, b}(z)  P(z)$,  we have:

\vspace{2mm}

\begin{equation}
H_{\sigma,f}(z)=\dfrac{P(z)}{S^0(z)}= \frac{{  \sqrt[4]{ (2 \pi  \sigma^2)} \big( a_{i,0}-(a_{i,1}  +a_{i,2}  +a_{i,3}) \big)}}{   \big( a_{i,0}-a_{i,1} z^{-1} -a_{i,2} z^{-2}-a_{i,3}z^{-3} \big)   }   \label {avanti3}
\end{equation}

\begin{equation}
H_{\sigma,b}(z)=\dfrac{S(z)}{P(z)}= \frac{ \sqrt[4]{ (2 \pi  \sigma^2)} {\big( a_{i,0}-(a_{i,1}  +a_{i,2}  +a_{i,3}) \big)}}{   \big( a_{i,0}-a_{i,1} z -a_{i,2} z^{2}-a_{i,3}z^{3} \big)   }   \label {avanti4} \qquad  \ \ \
\end{equation}

\vspace{2mm}

\noindent from which we obtain respectively:  

\vspace{2mm}
\footnotesize
\begin{eqnarray}
 \big( a_{i,0}  -a_{i,1} z^{-1} -a_{i,2} z^{-2} -a_{i,3}z^{-3}\big) P(z)    =\sqrt[4]{ (2 \pi  \sigma^2)} \big(a_{i,0}-(a_{i,1}  +a_{i,2}  +a_{i,3}) \big) S^0(z), 
\label{dd1}\\[3mm]
 \big (a_{i,0}  -a_{i,1} z^{-1}   -a_{i,2} z^{-2}  -a_{i,3}z^{-3}\big) S(z)    =\sqrt[4]{ (2 \pi  \sigma^2)} \big( a_{i,0}-(a_{i,1}  +a_{i,2}  +a_{i,3}) \big) P(z). \quad  
\label{dd2}
\end{eqnarray}
\normalsize

\vspace{2mm}

\noindent Antitransforming, by means of the  $\mathcal{Z}^{-1}$  transformation, the equations  (\ref{dd1}) and (\ref{dd2})   and by  the theorem of the delay  \citep[e.g.,][]{Oppe}, we obtain the following forward and backward  finite difference equations:

\vspace{2mm}

\begin{eqnarray}
p_i=\beta s^0_i+\alpha_{i,1} p_{i-1}+\alpha_{i,2} p_{i-2} +\alpha_{i,3} p_{i-3}.
\label{forward}\\[3mm]
s_i=\beta p_i+\alpha_{i,1} s_{i+1}+\alpha_{i,2} s_{i+2} +\alpha_{i,3} s_{i+3}.\ 
\label{backward}
\end{eqnarray}
\vspace{2mm}
 
\noindent where  the functions $p_i$ and $s_i$ are respectively the  $\mathcal{Z}^{-1}$  transformations of $P (z)$ and $S(z)$ functions and the filtering  coefficients  are  $\alpha_{i,1}=a_{i,1}/a_{i,0} $,  $\alpha_{i,2}=a_{i,2}/a_{i,0} $,  $\alpha_{i,3}=a_{i,3}/a_{i,0} $  and $\beta_i = \sqrt[4]{ (2 \pi  \sigma^2)}\big(1-(\alpha_{i,1}  +\alpha_{i,2}  +\alpha_{i,3})\big)$. Remembering that $\sigma= R_i/\delta x_i$ and observing the equations (\ref{angela}), (\ref{forward}) and (\ref{backward}) and  the values of the filtering coefficients $\alpha_{i,1}, \alpha_{i,2},  \alpha_{i,3}$ and $\beta_i$, then we obtain the thesis. 
$ $ \hfill $\blacksquare$
\end{proof}
}
\vspace{2mm}

\begin{remark}[Remark 1]
{\em
For an arbitrary input  distribution $s^0_i$,  
a measure for   accuracy   of  a   RF is given by the following inequality:  
\begin{equation} 
\|            \epsilon_{s_i}  \|_2 \leq \|           \epsilon_{h_i} \|_2    \|          s^0_i  \|_2 
\label{ineq}
\end{equation}
where:
\begin{itemize}
\item  {$\|            \epsilon_{s_i} \|_2$ is  the euclidean norm of  the difference 
between the discrete  convolution $s^*_i=g_i  \otimes s^0_i  $(with $g_i$ normalized  gaussian function)  and the function   $s_i$, obtained by the   RF applied to $s^0_i$.}
\item   { $\|            \epsilon_{h_i} \|_2$   is    the euclidean norm of  the difference 
between the  gaussian  $g_i$ and   the function  $h_i$ (called \emph{impulse response}),   obtained by a RF applied to  Dirac  function;} 
\item   { $\| s^0_i\|_2$ is the euclidean norm  of the input function $s^0_i$ .  }
\end{itemize}
}

\end{remark}

{\footnotesize
\begin{proof}
 The error in the output  per sample  when  we substitute  the true convolution $s^*_i$ with an  approximation $s_i$   is given by
\begin{equation}
{\epsilon}_{s_i} =          s^*_i-   s_i   =(           g_i -   h_i ) \otimes            s^0_i =             \varepsilon_{h_i}  \otimes            s^0_i . 
\end{equation}
     \noindent  Hence  it holds the following: 
\begin{equation}
\|            \epsilon _{s_i} \|_2 = \|           \epsilon _{h_i}  \otimes            s^0_i  \|_2   
\label{eps}
\end{equation}
\noindent and applying the Cauchy-Schwarz inequality   to the right-hand side of equation  (\ref{eps}), it  gives the thesis: 
\begin{equation}
\|            \epsilon_{s_i} \|_2 \leq \|           \epsilon_{h_i} \|_2    \|          s_i  \|_2   
\end{equation}
$ $ \hfill $\blacksquare$
\end{proof}
}

\noindent \textbf{References}

\bibliographystyle{elsarticle-harv}


\end{document}

\newpage
\newpage
\bibliographystyle{elsarticle-harv}






%
%
}
\end{document}